\documentclass{article}
\usepackage{graphicx}
\usepackage{amssymb}
\usepackage{amsfonts}
\usepackage{amsmath}

\begin{document}

\title{Finite Element Method for Cosserat Plates}
\author{Roman Kvasov$^{\ast }$ and Lev Steinberg$^{\ast \ast }$ \\
\\
$^{\ast }$Department of Mathematics\\
University of Puerto Rico at Aguadilla \\
Aguadilla, Puerto Rico 00604, USA\\ \\
$^{\ast \ast}$Department of Mathematical Sciences\\ 
University of Puerto Rico at Mayag\"{u}ez,\\
Mayag\"{u}ez, Puerto Rico 00681, USA}

\maketitle

\begin{abstract}
This paper presents the Finite Element Method for Cosserat plates. The mathematical model for Cosserat elastic plates is based on the calculation of the optimal value of the splitting parameter. We discuss the existence and uniqueness of the weak solution and the convergence of the proposed FEM. The Finite Element analysis of the clamped Cosserat plates of different shapes under different loads is provided. We present the numerical validation of the proposed FEM by estimating the order of convergence, when comparing the main kinematic variables with the analytical solution. We also consider the numerical analysis of plates with circular holes. We show that as expected the stress concentration factor around the hole is smaller than the classical value and smaller holes exhibit less stress concentration compared to larger ones.

\end{abstract}

\vspace{.25in}

\section{Introduction}

\vspace{.25in}

A complete theory of asymmetric elasticity introduced by the Cosserat brothers \cite{Cosserat1909} gave rise to a varierty of beam, shell and plate theories. The first theories of plates that take into account the microstructure of the material were developed in the 1960s. Eringen proposed a complete theory of plates in the framework of Cosserat (micropolar) Elasticity \cite{Eringen1967}, while independently Green and Naghdi specialized their general theory of Cosserat surface to obtain the linear Cosserat plate \cite{Green1966}. Numerous plate theories were formulated afterwards; for the extensive review of the latest developments we recommend to turn to \cite{Altenbach2010}.

The first theory of Cosserat elastic plates based on the Reissner plate theory was developed in \cite{Steinberg2010} and its finite element modeling is provided in \cite{Kvasov2011}. The enhanced version of the Cosserat plate theory was presented by the authors in \cite{Steinberg2013} and includes additional assumptions leading to the introduction of the splitting parameter. The theory provides the equilibrium equations and constitutive relations and the optimal value of the minimization of the elastic energy of the Cosserat plate. The paper also provides the analytical solutions of the presented plate theory and the three-dimensional Cosserat Elasticity for simply supported rectangular plate. The comparison of these solutions showed that the precision of the developed Cosserat plate theory is compatible with the precision of the Reissner plate theory.

The numerical modeling of bending of simply supported rectangular plates is given in \cite{Kvasov2013}. The paper provides the Cosserat plate field equations and the rigorous formula for the optimal value of the splitting parameter. The solution of the Cosserat plate converges to the Reissner plate theory \cite{Reissner1944}, \cite{Reissner1945} as the elastic asymmetric parameters tend to zero. The Cosserat plate theory shows agreement with the size-effect, confirming that the plates of smaller thickness are more rigid than expected from the Reissner model. The modeling of Cosserat plates with simply supported rectangular holes is also provided.

The extension of the static model of Cosserat elastic plates to dynamic problems is presented in \cite{Steinberg2015}. The computations predict a new kind of natural frequencies associated with the material microstructure and were shown to be consistent with the size-effect principle known from the Cosserat plate deformation reported in \cite{Kvasov2013}.

In this paper we present the Finite Element Method for Cosserat elastic plates based on the enhanced Cosserat plate theory given in \cite{Steinberg2013}. Since \cite{Kvasov2013} was restricted only to the case of rectangular plates, the current article represents an extension of this work for the Finite Element analysis of the Cosserat plates of different shapes, under different loads and different boundary conditions. We discuss the existence and uniqueness of the weak solution and the convergence of the proposed FEM. We present the numerical validation of the proposed FEM by estimating the order of convergence, when comparing the main kinematic variables with the analytical solution of the two-dimensional problem. We also consider the numerical analysis of plates with circular holes. We numerically calculate the stress concentration factor around the hole and show that it is smaller would be expected on the basis of Reissner theory for simple elastic plates. The finite element comparison of the plates with holes confirm that smaller holes exhibit less stress concentration.

\vspace{.25in}

\section{Cosserat Plate Equations}

\vspace{.25in}

In this section we will review the main equations of the Cosserat Plate Theory presented in \cite{Steinberg2013}. 

Throughout this article Greek indices are assumed to range from 1 to 2, while the Latin indices range from 1 to 3 if not specified otherwise. We will also employ the Einstein summation convention according to which summation is implied for any repeated index.

We will consider the thin plate $P$ of thickness $h$ and $x_{3}=0$ containing its middle plane. The sets $T$ and $B$ are the top and bottom surfaces contained in the planes $x_{3}=h/2$, $x_{3}=-h/2$ respectively and the curve $\Gamma$ is the boundary of the middle plane of the plate. The set of points $P=\left( \Gamma \times \lbrack -\frac{h}{2},\frac{h}{2}]\right) \cup T\cup B$ forms the entire surface of the plate. $\Gamma_{u}\times \lbrack -\frac{h}{2},\frac{h}{2}]$ is the lateral part of the boundary where displacements and microrotations are prescribed, while $\Gamma _{\sigma }\times \lbrack -\frac{h}{2},\frac{h}{2}]$ is the lateral part of the boundary edge where stress and couple stress are prescribed.

The assumptions on the displacements $u_{i}$ and microrotations $\phi{i}$ are given as
\begin{eqnarray}
u_{\alpha } &=& \zeta \Psi_{\alpha }(x_{1},x_{2}), \label{displacement_1} \\
u_{3} &=& W(x_{1},x_{2})+\left( 1-\zeta ^{2} \right) W^{\ast}(x_{1},x_{2}),
\label{displacement_2} \\
\phi _{\alpha } &=& \Omega_{\alpha }^{0}(x_{1},x_{2})\left(1-\zeta ^{2}\right) + \hat{\Omega}_{\alpha}(x_{1},x_{2}), 
\label{microrotation_1} \\
\phi_{3} &=& \zeta \Omega_{3}(x_{1},x_{2}), \label{microrotation_2}
\end{eqnarray}
where $\zeta =\frac{2x_{3}}{h}$ and $\alpha ,\beta \in \{1,2\}$.

The equilibrium system of equations for Cosserat plate bending is given as
\begin{eqnarray}
M_{\alpha \beta ,\alpha }-Q_{\beta } &=&0,  \label{bending_system_1} \\
Q_{\alpha ,\alpha }^{\ast }+\hat{p}_{1} &=&0,  \label{bending_system_2} \\
R_{\alpha \beta ,\alpha }+\varepsilon _{3\beta \gamma }\left( Q_{\gamma
}^{\ast }-Q_{\gamma }\right) &=&0,  \label{bending_system_3} \\
\varepsilon _{3\beta \gamma }M_{\beta \gamma }+S_{\alpha ,\alpha }^{\ast }
&=&0,  \label{bending_system_4} \\
\hat{Q}_{\alpha ,\alpha }+\hat{p}_{2} &=&0,  \label{bending_system_5} \\
R_{\alpha \beta ,\alpha }^{\ast }+\varepsilon _{3\beta \gamma }\hat{Q}%
_{\gamma } &=&0,  \label{bending_system_6}
\end{eqnarray}
where $M_{11}$ and $M_{22}$ are the bending moments, $M_{12}$ and $M_{21}$ -- twisting moments, $Q_{\alpha}$ -- shear forces, $Q_{\alpha}^{\ast}$, $\hat{Q}_{\alpha }$ -- transverse shear forces, $R_{11}$, $R_{22}$, $R_{11}^{\ast}$, $R_{22}^{\ast}$ -- micropolar bending moments, $R_{12}$, $R_{21}$, $R_{12}^{\ast}$ ,$R_{21}^{\ast}$ -- micropolar twisting moments, $S_{\alpha}^{\ast}$ -- micropolar couple moments, all defined per unit length. The initial pressure $p$ is represented here by the pressures $\hat{p}_{1}=\eta p$ and $\hat{p}_{2}=\frac{2}{3}\left( 1-\eta \right) p$, where $\eta $ is the splitting parameter.

It was shown that the system of equilibrium equations is accompanied by the zero variation of the stress energy with respect to the splitting parameter
\begin{equation}
\delta U_{K}^{S}(\eta)=0. \label{zero_variation}
\end{equation}

The constitutive formulas for Cosserat plate given in the following reverse form \cite{Steinberg2013}:
\begin{eqnarray}
M_{\alpha \alpha } &=&\frac{h^{3}\mu (\lambda +\mu )}{3(\lambda +2\mu )}\Psi_{\alpha ,\alpha }+\frac{\lambda \mu h^{3}}{6(\lambda +2\mu )}\Psi _{\beta,\beta }+\frac{\left( 3p_{1}+5p_{2}\right) \lambda h^{2}}{30(\lambda +2\mu )},\label{constitutive_formulas_reverse_1} \\
M_{\beta \alpha } &=&\frac{\left( \mu -\alpha \right) h^{3}}{12}\Psi_{\alpha ,\beta }+\frac{h^{3}(\alpha +\mu )}{12}\Psi _{\beta ,\alpha}+(-1)^{\beta }\frac{\alpha h^{3}}{6}\Omega _{3},
\label{constitutive_formulas_reverse_2} \\
R_{\beta \alpha } &=&\frac{5\left( \gamma -\epsilon \right) h}{6}\Omega_{\beta ,\alpha }^{0}+\frac{5h\left( \gamma +\epsilon \right) }{6}\Omega_{\alpha ,\beta }^{0}, \label{constitutive_formulas_reverse_3} \\
R_{\alpha \alpha } &=&\frac{10h\gamma \left( \beta +\gamma \right) }{3\left(\beta +2\gamma \right) }\Omega _{\alpha ,\alpha }^{0}+\frac{5h\beta \gamma }{3(\beta +2\gamma )}\Omega _{\beta ,\beta }^{0},\label{constitutive_formulas_reverse_4} \\
R_{\beta \alpha }^{\ast } &=&\frac{2\left( \gamma -\epsilon \right) h}{3}\hat{\Omega}_{\beta ,\alpha }+\frac{2\left( \gamma +\epsilon \right) h}{3}\hat{\Omega}_{\alpha ,\beta },  \label{constitutive_formulas_reverse_5} \\
R_{\alpha \alpha }^{\ast } &=&\frac{8\gamma \left( \gamma +\beta \right) h}{3(\beta +2\gamma )}\hat{\Omega}_{\alpha ,\alpha }+\frac{4\gamma \beta h}{3(\beta +2\gamma )}\hat{\Omega}_{\beta ,\beta },\label{constitutive_formulas_reverse_6} \\
Q_{\alpha } &=&\frac{5h(\alpha +\mu )}{6}\Psi _{\alpha }+\frac{5\left( \mu-\alpha \right) h}{6}W_{,\alpha }+\frac{2\left( \mu -\alpha \right) h}{3}W_{,\alpha }^{\ast }  \notag \\
&&+(-1)^{\beta }\frac{5h\alpha }{3}\left( \Omega _{\beta }^{0}+\hat{\Omega}_{\beta }\right),\label{constitutive_formulas_reverse_7} \\
Q_{\alpha }^{\ast } &=&\frac{5\left( \mu -\alpha \right) h}{6}\Psi _{\alpha}+\frac{5\left( \mu -\alpha \right) ^{2}h}{6\left( \mu +\alpha \right) }W_{,\alpha }+\frac{2\left( \mu +\alpha \right) h}{3}W_{,\alpha }^{\ast } \notag \\
&&+(-1)^{\alpha }\frac{5h\alpha }{3}\left( \Omega _{\beta }^{0}+\frac{\left(\mu -\alpha \right) }{\left( \mu +\alpha \right) }\hat{\Omega}_{\beta}\right), \label{constitutive_formulas_reverse_8} \\
\hat{Q}_{\alpha } &=&\frac{8\alpha \mu h}{3\left( \mu +\alpha \right)}W_{,\alpha }+(-1)^{\alpha }\frac{8\alpha \mu h}{3\left( \mu +\alpha \right)}\hat{\Omega}_{\beta},  \label{constitutive_formulas_reverse_9} \\
S_{\alpha }^{\ast } &=&\frac{5\gamma \epsilon h^{3}}{3\left( \gamma
+\epsilon \right) }\Omega _{3,\alpha}.\label{constitutive_formulas_reverse_10}
\end{eqnarray}
In these formulas the greek subindex $\beta=1$ iff $\alpha=2$ and $\beta=2 $ iff $\alpha=1$. The parameters $\lambda$ and $\mu$ are the Lam\'e constants and $\alpha$, $\beta$, $\gamma$ and $\epsilon$ are asymmetric constants.

In order to obtain the micropolar plate bending field equations in terms of the kinematic variables, the constitutive formulas in the reverse form (\ref{constitutive_formulas_reverse_1}) - (\ref{constitutive_formulas_reverse_10}) are substituted into the bending system of equations (\ref{bending_system_1}) - (\ref{bending_system_6}). The obtained Cosserat plate bending field equations can be represented as an elliptic system of nine partial differential equations in terms of the kinematic variables \cite{Kvasov2013}:
\begin{equation}
L v = f \left( \eta \right) \label{field_equations}
\end{equation}
where $L$ is a linear differential operator acting on the vector of  kinematic variables $v$ (unknowns), and $f\left( \eta \right)$ is the right-hand side vector defined as (\ref{righthand_side_vector}), that in general depends on $\eta$:
\begin{equation}
L=\left[ 
\begin{array}{ccccccccc}
L_{11} & L_{12} & L_{13} & L_{14} & 0 & L_{16} & k_{1}L_{13} & 0 & L_{16} \\ 
L_{12} & L_{22} & L_{23} & L_{24} & L_{16} & 0 & k_{1}L_{23} & L_{16} & 0 \\ 
-L_{13} & -L_{23} & L_{33} & 0 & L_{35} & L_{36} & k_{1}L_{77} & L_{38} & L_{39} \\ 
L_{41} & L_{42} & 0 & L_{44} & 0 & 0 & 0 & 0 & 0 \\ 
0 & -L_{16} & -L_{38} & 0 & L_{55} & L_{56} & -k_{1}L_{35} & L_{58} & 0 \\ 
L_{16} & 0 & -L_{39} & 0 & L_{56} & L_{66} & -k_{1}L_{36} & 0 & L_{58} \\ 
-L_{13} & -L_{14} & L_{73} & 0 & L_{35} & L_{36} & k_{1}L_{77} & L_{78} & L_{79} \\ 
0 & -L_{16} & -L_{78} & 0 & L_{85} & L_{56} & -k_{1}L_{35} & k_{1}L_{88} & k_{1}L_{56} \\ 
L_{16} & 0 & -L_{79} & 0 & L_{56} & L_{55} & -k_{1}L_{36} & k_{1}L_{56} & k_{1}L_{99}
\end{array}
\right],\label{differential_operator}
\end{equation}

\begin{equation}
v =\left[ 
\begin{array}{ccccccccc}
\Psi _{1}, & \Psi _{2}, & W, & \Omega_{3}, & \Omega_{1}^{0}, & \Omega_{2}^{0}, & W^{\ast }, & \Omega_{1}^{0}, & \Omega_{2}^{0}%\
\end{array}%
\right] ^{T},\label{solution_vector}
\end{equation}

\begin{equation}
f\left( \eta \right) =\left[ 
\begin{array}{ccccccccc}
-\frac{h^{2}\lambda \left( 3p_{1,1}+5p_{2,1}\right) }{30\left( \lambda+2\mu \right) }, & -\frac{h^{2}\lambda \left(3p_{1,2}+5p_{2,2}\right) }{30\left( \lambda +2\mu \right) }, & -p_{1}, & 0, & 0, & 0, & \frac{h^{2}(3p_{1}+4p_{2})}{24}, & 0, & 0
\end{array}\right]^{T}\label{righthand_side_vector}
\end{equation}

The operators $L_{ij}$ are defined as follows

$L_{11}=c_{1}\frac{\partial ^{2}}{\partial x_{1}^{2}}+c_{2}\frac{\partial^{2}}{\partial x_{2}^{2}}-c_{3},L_{12}=(c_{1}-c_{2})\frac{\partial ^{2}}{\partial x_{1}x_{2}},L_{13}=c_{11}\frac{\partial }{\partial x_{1}},L_{14}=c_{12}\frac{\partial }{\partial x_{2}},L_{16}=c_{13},L_{17}=k_{1}c_{11}\frac{\partial }{\partial x_{1}},L_{22}=c_{2}\frac{\partial ^{2}}{\partial x_{1}^{2}}+c_{1}\frac{\partial ^{2}}{\partial x_{2}^{2}}-c_{3},L_{23}=c_{11}\frac{\partial }{\partial x_{2}},L_{24}=-c_{12}\frac{\partial }{\partial
x_{1}},L_{33}=c_{3}(\frac{\partial ^{2}}{\partial x_{1}^{2}}+\frac{\partial
^{2}}{\partial x_{2}^{2}}),L_{35}=-c_{13}\frac{\partial }{\partial x_{2}}%
,L_{36}=c_{13}\frac{\partial }{\partial x_{1}},L_{38}=-c_{10}\frac{\partial 
}{\partial x_{2}},L_{39}=c_{10}\frac{\partial }{\partial x_{1}}%
,L_{41}=-c_{12}\frac{\partial }{\partial x_{2}},L_{42}=c_{12}\frac{\partial 
}{\partial x_{1}},L_{44}=c_{6}\left( \frac{\partial ^{2}}{\partial x_{1}^{2}}%
+\frac{\partial ^{2}}{\partial x_{2}^{2}}\right) -2c_{12},L_{55}=c_{7}\frac{%
\partial ^{2}}{\partial x_{1}^{2}}+c_{8}\frac{\partial ^{2}}{\partial
x_{2}^{2}}-2c_{13},L_{56}=(c_{7}-c_{8})\frac{\partial ^{2}}{\partial
x_{1}x_{2}},L_{58}=-c_{9},L_{66}=c_{8}\frac{\partial ^{2}}{\partial x_{1}^{2}}+c_{7}\frac{\partial
^{2}}{\partial x_{2}^{2}}-2c_{13},L_{73}=c_{5}(\frac{\partial ^{2}}{\partial
x_{1}^{2}}+\frac{\partial ^{2}}{\partial x_{2}^{2}}),L_{77}=c_{4}(%
\frac{\partial ^{2}}{\partial x_{1}^{2}}+\frac{\partial ^{2}}{\partial
x_{2}^{2}}),L_{78}=-c_{14}\frac{\partial }{\partial x_{2}},L_{79}=c_{14}%
\frac{\partial }{\partial x_{1}},L_{85}=c_{7}\frac{\partial ^{2}}{\partial
x_{1}^{2}}+c_{8}\frac{\partial ^{2}}{\partial x_{2}^{2}}%
-2c_{13},L_{88}=c_{7}\frac{\partial ^{2}}{\partial x_{1}^{2}}+c_{8}\frac{%
\partial ^{2}}{\partial x_{2}^{2}}-c_{15},L_{99}=c_{8}\frac{\partial ^{2}}{%
\partial x_{1}^{2}}+c_{7}\frac{\partial ^{2}}{\partial x_{2}^{2}}-c_{15}$

The coefficients $c_{i}$ are given as

$c_{1}=\frac{h^{3}\mu (\lambda +\mu )}{3(\lambda +2\mu )},c_{2}=\frac{h^{3}(\alpha +\mu )}{12},c_{3}=\frac{5h(\alpha +\mu )}{6},c_{4}=\frac{5h(\alpha-\mu)^{2}}{6(\alpha+\mu )},c_{5}=\frac{h(5\alpha ^{2}+6\alpha \mu +5\mu ^{2})}{6(\alpha +\mu )}
,c_{6}=\frac{h^{3}\gamma \epsilon }{3(\gamma +\epsilon )},c_{7}=\frac{10h\gamma \left( \beta +\gamma \right) }{3\left( \beta +2\gamma \right)},c_{8}=\frac{5h\left( \gamma +\epsilon \right) }{6}c_{9}=\frac{10h\alpha^{2}}{3(\alpha +\mu )},c_{10}=\frac{5h\alpha (\alpha -\mu )}{3(\alpha +\mu )},c_{11}=\frac{5h(\alpha -\mu )}{6},c_{12}=\frac{h^{3}\alpha }{6},c_{13}=\frac{5h\alpha }{3},c_{14}=\frac{h\alpha (5\alpha +3\mu )}{3(\alpha +\mu )},c_{15}=\frac{2h\alpha (5\alpha+4\mu )}{3(\alpha +\mu )}$.

The optimal value of the splitting parameter $\eta$ is given as
\begin{equation}
\eta _{0}=\frac{2\mathcal{W}^{\left( 00\right) }-\mathcal{W}^{\left(
10\right) }-\mathcal{W}^{\left( 01\right) }}{2\left( \mathcal{W}^{\left(
11\right) }+\mathcal{W}^{\left( 00\right) }-\mathcal{W}^{\left( 10\right) }-%
\mathcal{W}^{\left( 01\right) }\right) },  \label{eta_minimum}
\end{equation}
where  $\mathcal{W}^{\left( ij\right)}$ are the work densities provided in \cite{Kvasov2013}.

\vspace{.25in}

\section{Finite Element Algorithm for Cosserat Plate}

\vspace{.25in}

The right-hand side of the system (\ref{field_equations}) depends on the splitting parameter $\eta$ and so does the solution $v$, that we will formally denote as $v_{\eta}$. Therefore the solution of the Cosserat elastic plate bending problem requires not only solving the system (\ref{field_equations}), but also an additonal technique for the calculation of the value of the splitting parameter, that corresponds to the unique solution. Considering that the elliptic systems of partial differential equations correspond to a state where the minimum of the energy is reached, the optimal value of the splitting parameter should minimize the elastic plate energy \cite{Solin2006}. The minimization corresponds to the zero variation of the plate stress energy (\ref{zero_variation}).

The Finite Element Method for Cosserat elastic plates is based on the algorithm for the optimal value of the splitting parameter. This algorithm requires solving the system (\ref{field_equations}) for two different values of the splitting parameter $\eta$, numerical calculation of stresses, strains and the corresponding work densities. We will follow \cite{Kvasov2013} in the description of our Finite Element Method algorithm:

1. Use classic Galerkin FEM to solve two elliptic systems: 
\begin{eqnarray}
L v_{0} = f\left( 0 \right) \notag \\
L v_{1} = f\left( 1 \right) \notag
\end{eqnarray}
for $v_{0}$ and $v_{1}$ respectively.

2. Calculate the optimal value of the splitting parameter $\eta _{0}$ using (\ref{eta_minimum}).

3. Calculate the optimal solution $v_{\eta _{0}}$ of the Cosserat plate bending problem as a linear combination of $v_{0}$ and $v_{1}$:
\begin{equation}
v_{\eta_{0}} = (1-\eta_{0})v_{0}+\eta_{0} v_{1}.\label{solution_linear_combination}
\end{equation}

\vspace{.25in}

\subsection{Weak Formulation of the Clamped Cosserat Plate}

\vspace{.25in}

Let us consider the following hard clamped boundary conditions similar to \cite{Arnold1989}:
\begin{eqnarray}
&&W=0, W^{\ast }=0, \mathbf{\Psi }\cdot \hat{\normalfont\textbf{s}} = 0, \mathbf{\Psi }\cdot \hat{\normalfont\textbf{n}} = 0, \Omega_{3} =0,  
\label{clamped_bc_1} \\
&& \mathbf{\Omega^{0}}\cdot \hat{\normalfont\textbf{s}} = 0, \mathbf{\Omega^{0}}\cdot \hat{\normalfont\textbf{n}} = 0, \mathbf{\hat{\Omega}}\cdot \hat{\normalfont\textbf{s}} = 0, \mathbf{\hat{\Omega}}\cdot \hat{\normalfont\textbf{n}} = 0,
\label{clamped_bc_2}
\end{eqnarray}
where $\hat{\normalfont\textbf{n}}$ and $\hat{\normalfont\textbf{s}}$ are the normal and the tangent vectors to the boundary. These conditions represent homogeneous Dirichlet type boundary conditions for the kinematic variables:
\begin{eqnarray}
&& W=0, \text{ }W^{\ast }=0, \text{ }\Psi _{1}=0, \text{ }\Psi _{2}=0, \text{ }\Omega_{3}=0, 
\label{clamped_boundary_1} \\
&& \hat{\Omega}_{1}^{0}=0, \text{ }\hat{\Omega}_{2}^{0}=0, \text{ }\hat{\Omega}_{1}^{0}=0, \text{ }\hat{\Omega}_{2}^{0}=0.
\label{clamped_boundary_2}
\end{eqnarray}

Let us denote by $\mathbf{L}^{2}\left(B_{0}\right)$ the standard space of square-integrable functions defined everywhere on $B_{0}$: 
\begin{equation*}
\mathbf{L}^{2}\left(B_{0}\right) = \left\{ v : \int_{B_{0}}{v^{2}ds} < \infty \right\}
\end{equation*}
and by $\mathbf{H}^{1}\left( B_{0}\right )$ the Hilbert space of functions that are square-integrable together with their first partial derivatives:
\begin{equation*}
\mathbf{H}^{1}\left( B_{0}\right) = \left\{v : v\in L^{2}\left(B_{0}\right), \partial_{i}v\in L^{2}\left(B_{0}\right) \right\}
\end{equation*}

Let us denote the Hilbert space of functions from $\mathbf{H}^{1}\left( B_{0}\right)$ that vanish on the boundary as in \cite{Johnson1987}:
\begin{equation*}
\mathbf{H}^{1}_{0}\left( B_{0} \right) = \left\{v\in \mathbf{H}^{1}\left( B_{0}\right) ,v=0 \text{ on } \partial B_{0} \right\} \end{equation*}

The space $\mathbf{H}^{1}_{0}\left( B_{0}\right)$ is equipped with the inner product:
\begin{equation*}
\left\langle u,v \right\rangle _{\mathbf{H}^{1}_{0}} = \int_{B_{0}}{\left(uv+\partial_{i}u\partial_{i}v\right)ds}
\text{ for } u,v \in\mathbf{H}^{1}_{0}\left(B_{0}\right)
\end{equation*}

Taking into account that the boundary conditions for all variables are of the same homogeneous Dirichlet type, we look for the solution in the function space $\mathcal{H}\left( B_{0}\right)$ defined as
\begin{equation}
\mathcal{H} = \mathbf{H}^{1}_{0} \left( B_{0}\right)^{9}.
\label{space_H}
\end{equation}

The space $\mathcal{H}$ is equipped with the inner product $\left\langle u,v \right\rangle_{\mathcal{H}}$:
\begin{equation*}
\left\langle u,v \right\rangle_{\mathcal{H}} = \sum_{i=1}^{9}{\left\langle u_{i},v_{i} \right\rangle_{\mathbf{H}^{1}_{0}}} \text{ for } u,v \in \mathcal{H}
\end{equation*}
and relative to the metric 
\begin{equation*}
d\left(u,v\right) = \left\| u-v \right\|_{\mathcal{H}} \text{ for } u,v \in \mathcal{H},
\end{equation*}
induced by the norm $\left\| x \right\| = \sqrt {\left\langle x,x \right\rangle_{\mathcal{H}}}$, the space $\mathcal{H}$ is a complete metric space and therefore is a Hilbert space \cite{Conway1985}.

Let us consider a dot product of both sides of the system of the field equations (\ref{field_equations}) and an arbitrary function $v \in \mathcal{H}$:
\begin{equation*}
v \cdot L u = v \cdot f \left( \eta \right)
\end{equation*}
and then integrate both sides of the obtained scalar equation over the plate $B_{0}$:
\begin{equation*}
\int_{B_{0}}\left(v \cdot L u \right)ds = \int_{B_{0}}\left(v \cdot f\left( \eta \right) \right)ds.
\end{equation*}

Let us introduce a bilinear form $a\left( u,v\right) : \mathcal{H} \times \mathcal{H} \rightarrow \mathbb{R}$ and a linear form $b_{(\eta)} \left( v\right) : \mathcal{H} \rightarrow \mathbb{R}$ defined as
\begin{eqnarray}
a\left( u,v\right) &=& \int_{B_{0}} \left(v \cdot L u \right)ds, \label{bilinear_form} \\
b_{(\eta)} \left( v\right) &=& \int_{B_{0}}\left(v \cdot f\left( \eta \right) \right)ds. \notag
\end{eqnarray}

The expression for $a\left( u,v\right)$
\begin{equation*}
a\left( v, u \right) = \int_{B_{0}}\left(v_{i}L_{ij}u_{j}\right)ds
\end{equation*}
is a summation over the terms of the form
\begin{equation*}
a^{ij} \left( v_{m}, u_{n} \right) = \int_{B_{0}}\left( v_{m} \hat{L}  u_{n} \right)ds,
\end{equation*}
where $v_{m} \in \mathcal{H}_{m}$, $u_{n} \in \mathcal{H}_{n}$ and $\hat{L}$ is a scalar differential operator. 

There are 3 types of linear operators present in the field equations (\ref{field_equations}) -- operators of order zero, one and two, which are constant multiples of the following differential operators:
\begin{eqnarray}
L^{\left( 0\right) } &=& {1},  \label{order_0_operator} \\
L^{\left( 1\right) } &=&\frac{\partial }{\partial x_{\alpha}}, \label{order_1_operator} \\
L^{\left( 2\right) } &=&- \nabla \cdot A\nabla, \label{order_2_operator}
\end{eqnarray}

These operators act on the components of the vector $u$ and are multiplied by the components of the vector $v$ and the obtained expressions are then integrated over $B_{0}$:
\begin{eqnarray}
\int_{B_{0}}\left(v_{m} L^{\left( 0\right) }u_{n}\right) ds &=&\int_{B_{0}}\left(v_{m} u_{n}\right) ds \label{weak_order_zero} \\
\int_{B_{0}}\left(v_{m} L^{\left( 1\right) }u_{n}\right) ds &=&\int_{B_{0}}\left(v_{m} \frac{\partial u_{n}}{\partial x_{\alpha}}\right)ds 
\label{weak_order_one} \\
\int_{B_{0}}\left(v_{m} L^{\left( 2\right) }u_{n}\right) ds &=&-\int_{B_{0}}\left(v_{m} ( \nabla \cdot A\nabla ) u_{n}\right) ds \notag
\end{eqnarray}
where $v_{m} \in \mathcal{H}_{m}$ and $u_{n} \in \mathcal{H}_{n}$.
 
The weak form of the second order operator is obtained by performing the corresponding integration by parts and taking into account that the test functions $v_{m}$ vanish on the boundary $\partial B_{0}$:
\begin{eqnarray}
\int_{B_{0}}\left(v_{m} L^{\left( 2\right) }u_{n}\right) ds &=&-\int_{B_{0}}\left(v_{m} ( \nabla \cdot A\nabla u_{n}) \right) ds \notag \\
&=& -\int_{\partial B_{0}}\left( A\nabla u_{n}\cdot n\right) v_{m} d\tau+\int_{B_{0}}\left(A\nabla u_{n}\cdot \nabla v_{m}\right) ds \notag \\
&=& \int_{B_{0}}\left(A\nabla u_{n}\cdot \nabla v_{m}\right) ds\label{weak_order_two}
\end{eqnarray}

The expression for $b_{(\eta)}\left( v\right)$:
\begin{equation*}
b_{(\eta)}\left( v\right) =  \int_{B_{0}} v_{i} f_{i} \left( \eta \right)ds
\end{equation*}
represents a summation over the terms of the form:
\begin{equation*}
\int_{B_{0}} v_{m} \hat{f} \left( \eta \right)ds,
\end{equation*}

Taking into account that the optimal solution of the field equations (\ref{field_equations}) minimizes the plate stress energy, we can give the weak formulation for the clamped Cosserat plate bending problem.

\vspace{.25in}

\textbf{Weak Formulation of the Clamped Cosserat Plate Bending Problem}

\bigskip

\textit{Find all $u\in \mathcal{H}$ and $\eta\in\mathbb{R}$ that minimize the plate stress energy $U_{K}^{\mathcal{S}}\left(u,\eta\right)$ subject to}
\begin{equation}
a\left( v,u \right) = b_{(\eta)}\left( v\right) \textit{ for all $v \in \mathcal{H}$}
\label{weak_formulation}
\end{equation}

\vspace{.25in}

\subsection{Construction of the Finite Element Spaces}

\vspace{.25in}

Let us construct the finite element space, i.e. finite-dimensional subspace $\mathcal{H}_{h}$ of the space $\mathcal{H}$, where we will be looking for an approximate Finite Element solution of the weak formulation (\ref{weak_formulation}). 

\begin{figure}
\begin{center}
	\includegraphics{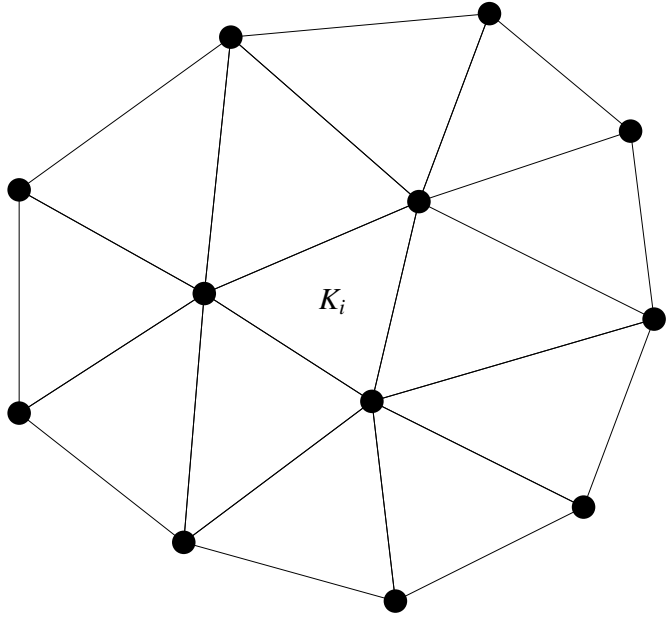}
	\vspace{.1in}
	\caption{Example of the Finite Element triangulation of the domain $B_{0}$}
	\label{fig:triangulation}
\end{center}
\end{figure}

Let us assume that the boundary $\partial B_{0}$ is a polygonal curve. Let us make a triangulation of the domain $B_{0}$ by subdividing $B_{0}$ into $l$ non-overlapping triangles $K_{i}$ with $m$ vertices $N_{j}$: 
\begin{equation*}
B_{0} = \bigcup_{i=1}^{l}{K_{i}}=K_{1}\cup K_{2}\cup ... \cup K_{l}
\end{equation*}
such that no vertex of the triangular element lies on the edge of another triangle (see Figure \ref{fig:triangulation}).
\noindent

Let us introduce the mesh parameter $h$ as the greatest diameter among the elements $K_{i}$:
\begin{equation*}
h = \max_{i=\overline{1,l}} d\left(K_{i}\right),
\end{equation*}
which for the triangular elements corresponds to the length of the longest side of the triangle.

We now define the finite dimensional space $\mathcal{\hat{H}}_{h}$ as a space of all continuous functions that are linear on each element $K_{j}$ and vanish on the boundary:
\begin{equation*}
\mathcal{\hat{H}}_{h} = \mathcal{H}_{i}^{h} = \left\{ v: v \in C\left(B_{0}\right), v \text{ is linear on every } K_{j}, v=0 \text{ on } \partial B_{0} \right\}.
\end{equation*}
By definition $\mathcal{H}_{i}^{h} \subset \mathcal{H}_{i}$, and the finite element space $\mathcal{H}_{h}$ is then defined as:
\begin{equation}
\mathcal{H}_{h} = \mathcal{\hat{H}}_{h}^{9}
\label{finite_space_H}
\end{equation}

The approximate weak solution $u^{h}$ can be found from the Galerkin formulation of the clamped Cosserat plate bending problem \cite{Hughes2004}.

\vspace{.25in}

\textbf{Galerkin Formulation of the Clamped Cosserat Plate}

\vspace{.25in}

\textit{Find all $u^{h}\in \mathcal{H}_{h}$ and $\eta\in\mathbb{R}$ that minimize the stress plate energy $U_{K}^{\mathcal{S}}\left(u^{h},\eta\right)$ subject to}
\begin{equation}
a\left( v^{h},u^{h} \right) = b_{(\eta)}\left( v^{h} \right) \textit{ for all $v^{h} \in \mathcal{H}_{h}$}
\label{galerkin_weak_formulation}
\end{equation}

\vspace{.25in}

The description of the function $v^{h}_{i} \in \mathcal{H}_{i}^{h}$ is provided by the values $v^{h}_{i} \left(N_{k}\right)$ at the nodes $N_{k}$ ($k=\overline{1,m}$).

Let us define the set of basis functions $\left\{ \phi_{1},\phi_{2},...,\phi_{m} \right\}$ of each space $\mathcal{H}_{i}^{h}$ as
\begin{equation*}
\phi_{j} \left(N_{k}\right) = \delta_{jk}, \text{ } j,k=\overline{1,m} 
\end{equation*}
excluding the points $N_{k}$ on the boundary $\partial B_{0}$.

\begin{figure}
\begin{center}
	\includegraphics{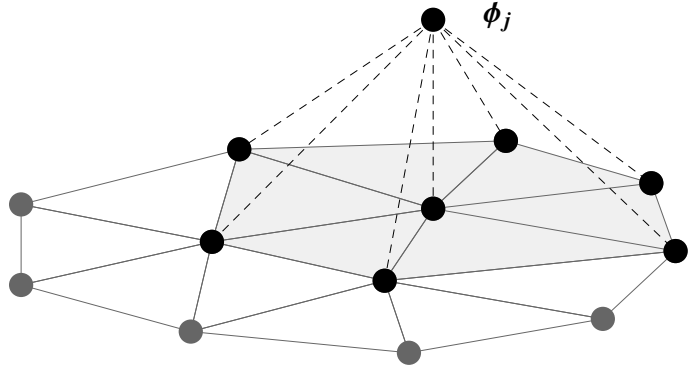}
	\vspace{.1in}
	\caption{Example of the Finite Element basis function}
	\label{fig:basis_function}
\end{center}
\end{figure}

Therefore
\begin{equation*}
\mathcal{H}_{i}^{h} = span \left\{\phi_{1},\phi_{2},...,\phi_{m}\right\} = \left\{ v : v = \sum_{j=1}^{m} {\alpha_{j}^{(i)} \phi_{j} } \right\}
\end{equation*}
and the functions $\phi_{j}$ is non-zero only at the node $N_{j}$ and those that belong to the specified boundary and the support of $\phi_{j}$ consists of all triangles $K_{i}$ with the common node $N_{j}$ (see the Figure \ref{fig:basis_function}).
\noindent

\vspace{.25in}

Since the spaces $\mathcal{H}_{i}^{h}$ are identical they will also have identical sets of basis functions $\phi_{j}$  ($j=\overline{1,m}$). Sometimes we will need to distinguish between the basis functions of different spaces assigning the superscript of the functions space to the basis function, i.e. the basis functions for the space $\mathcal{H}_{i}^{h}$ are $\phi_{j}^{i}$. For computational purposes these superscripts will be droped.

\vspace{.25in}

\subsection{Calculation of the Stiffness Matrix and the Load Vector}

\vspace{.25in}

The bilinear form of the Galerkin formulation (\ref{galerkin_weak_formulation}) is given as 
\begin{equation}
a\left( v^{h}, u^{h} \right) = a^{ij} \left( v_{i}^{h}, u_{j}^{h} \right) = \int_{B_{0}} v_{i}^{h} L_{ij} u_{j}^{h} ds
\label{summation_of_weak_forms}
\end{equation}

Since $u^{h}_{j} \in \mathcal{H}_{j}^{h}$ then there exist such constants $\alpha_{p}^{(j)} \in \mathbb{R}$ that
\begin{equation*}
u^{h}_{j} = \alpha_{p}^{(j)} \phi^{(i)}_{p}
\end{equation*}

Since the equation (\ref{summation_of_weak_forms}) is satisfied for all $v_{i}^{h} \in \mathcal{H}_{i}^{h}$ then it is also satisfied for all basis functions $\phi^{(i)}_{k}$ ($k=\overline{1,m}$):
\begin{equation*}
a^{ij} \left( v_{i}^{h}, u_{j}^{h} \right) = a^{ij} \left( \phi^{(i)}_{k}, \alpha_{p}^{(j)} \phi^{(j)}_{p} \right) = \alpha_{p}^{(j)} a^{ij} \left( \phi_{k}^{(i)}, \phi_{p}^{(j)} \right)
\end{equation*}
where 
\begin{equation}
a^{ij} \left( v, u \right) = \int_{B_{0}} v L_{ij} u ds
\label{block_weak_forms}
\end{equation}

Following \cite{Hughes2} we define the block stiffness matrices $K^{ij}$ ($i,j=\overline{1,9}$):
\[
K^{ij}=\left[ 
\begin{array}{ccc}
a^{ij}\left( \phi _{1}^{(i)},\phi _{1}^{(j)}\right)  & \ldots  & a^{ij}\left( \phi
_{1}^{(i)},\phi _{m}^{(j)}\right)  \\ 
\vdots  & \ddots  & \vdots  \\ 
a^{ij}\left( \phi _{m}^{(i)},\phi _{1}^{(j)}\right)  & \ldots  & a^{ij}\left( \phi
_{m}^{(i)},\phi _{m}^{(j)}\right) 
\end{array}%
\right]
\]

For computational purposes the superscripts of the basis functions can be droped and the block stiffness matrices $K^{ij}$ can be calculated as
\[
K^{ij}=\left[ 
\begin{array}{ccc}
a^{ij}\left( \phi _{1},\phi _{1}\right)  & \ldots  & a^{ij}\left( \phi
_{1},\phi _{m}\right)  \\ 
\vdots  & \ddots  & \vdots  \\ 
a^{ij}\left( \phi _{m},\phi _{1}\right)  & \ldots  & a^{ij}\left( \phi
_{m},\phi _{m}\right) 
\end{array}%
\right] 
\]

Let us define the block load vectors $F^{i}(\eta)$ ($i=\overline{1,9}$):
\[
F^{i}(\eta)=\left[ 
\begin{array}{c}
b_{(\eta)}^{i}\left( \phi _{1}\right)  \\ 
\vdots  \\ 
b_{(\eta)}^{i}\left( \phi _{m}\right) 
\end{array}%
\right] 
\]
and the solution block vectors $\alpha^{i}$ corresponding to the variable $u_{i}^{h}$ ($i=\overline{1,9}$):
\[
\alpha^{i}=\left[ 
\begin{array}{c}
\alpha _{1}^{i} \\ 
\vdots  \\ 
\alpha _{m}^{i}%
\end{array}%
\right] 
\]

The equation (\ref{galerkin_weak_formulation}) of the Galerkin formulation can be rewritten as
\begin{equation}
\left(K^{ij}\right) \alpha^{i} = F^{j}(\eta)
\label{weak_block_formulation}
\end{equation}

The global stiffness matrix consists of 81 block stiffness matrices $K^{ij}$, the global load vector consists of 9 block load vectors $F^{i}(\eta)$ and the global displacement vector is represented by the 9 blocks of coefficients $\alpha^{i}$. The entries of the block matrices $K^{ij}$ and the block vectors $F^{i}(\eta)$ can be calculated as
\begin{eqnarray}
K_{mn}^{ij} &=& \int_{B_{0}} \phi_{m} L_{ij} \phi_{n}  ds \notag \\
F_{m}^{i} (\eta) &=& \int_{B_{0}} \phi_{m} f_{i} \left( \eta \right) ds \notag
\end{eqnarray}

The block matrix form of the equation (\ref{galerkin_weak_formulation}) is given as

$\left[ 
\begin{array}{cccc}
\boxed{
\begin{array}{ccc}
&  &  \\ 
& K^{11} &  \\ 
&  & 
\end{array}
}
& 
\boxed{
\begin{array}{ccc}
&  &  \\ 
& K^{12} &  \\ 
&  & 
\end{array}
}
& \ldots  & 
\boxed{
\begin{array}{ccc}
&  &  \\ 
& K^{19} &  \\ 
&  & 
\end{array}
}
\\ 
\boxed{
\begin{array}{ccc}
&  &  \\ 
& K^{21} &  \\ 
&  & 
\end{array}
}
& 
\boxed{
\begin{array}{ccc}
&  &  \\ 
& K^{22} &  \\ 
&  & 
\end{array}
}
& \ldots  & 
\boxed{
\begin{array}{ccc}
&  &  \\ 
& K^{29} &  \\ 
&  & 
\end{array}
}
\\ 
\vdots  & \vdots  & \ddots  & \vdots  \\ 
\boxed{
\begin{array}{ccc}
&  &  \\ 
& K^{91} &  \\ 
&  & 
\end{array}
}
& 
\boxed{
\begin{array}{ccc}
&  &  \\ 
& K^{92} &  \\ 
&  & 
\end{array}
}
& \ldots  & 
\boxed{
\begin{array}{ccc}
&  &  \\ 
& K^{99} &  \\ 
&  & 
\end{array}%
}
\end{array}%
\right] 
\left[ 
\begin{array}{c}
\boxed{
\begin{array}{c}
\\ 
\alpha^{1} \\
\text{ } 
\end{array}
}
\\ 
\boxed{
\begin{array}{c}
\\ 
\alpha^{2} \\
\text{ }  
\end{array}
}
\\ 
\vdots  \\ 
\boxed{
\begin{array}{c}
\\ 
\alpha^{9} \\
\text{ }  
\end{array}
}
\end{array}
\right] =
\left[ 
\begin{array}{c}
\boxed{
\begin{array}{c}
\\ 
F^{1}(\eta) \\
\text{ }  
\end{array}
}
\\ 
\boxed{
\begin{array}{c}
\\ 
F^{2}(\eta) \\
\text{ }  
\end{array}
}
\\ 
\vdots  \\ 
\boxed{
\begin{array}{c}
\\ 
F^{9}(\eta) \\
\text{ }  
\end{array}%
}
\end{array}%
\right] $

\vspace{.25in}

\subsection{Existence and Convergence Remarks}

\vspace{.25in}

We will follow \cite{Costabel2010} and \cite{Costabel2012}, where the analysis of the analytic regularity for the linear elliptic systems and their general treatment were recently presented.

Let us consider the bilinear form $a(\cdot,\cdot)$ defined in (\ref{bilinear_form}):
\begin{equation*}
a\left( u,v\right) = \int_{B_{0}} v_{i} L_{ij} u_{j} ds
\end{equation*}
where $L_{ij}$ are linear differential operators of at most second order. Employing integration by parts for the second order operators $L_{ij}$ the bilinear form $a(\cdot,\cdot)$ can be rewritten in the following form:
\begin{equation*}
a\left( u,v\right) = \sum_{\left| \beta \right|, \left| \gamma \right| \leq 1} {\int_{B_{0}} c_{ij} \partial^{\beta} v_{i} \partial^{\gamma} u_{j} ds}
\end{equation*}
where $\beta$ and $\gamma$ are multi-indices.

Since coefficients $c_{ij}$ are constant and therefore bounded on $B_{0}$, the bilinear form $a(\cdot,\cdot)$ is continuous over $\mathcal{H}$ \cite{Costabel2010}, i.e. there exists a constant $C>0$ such that 
\begin{equation*}
\left|a(v,u)\right| \leq C \left\| v \right\|_{\mathcal{H}} \left\| u \right\|_{\mathcal{H}} \hspace{0.25in} \forall u,v \in \mathcal{H}. 
\end{equation*}

The strong ellipticity of the operator $L$ was shown in\cite{Kvasov2013}. Since the operator $L$ is strong elliptic on $B_{0}$ the bilinear form $a(\cdot,\cdot)$ is V-elliptic on $\mathcal{H}$ \cite{Costabel2010}, \cite{Solin2006}, i.e. there exists a constant $\alpha>0$ such that
\begin{equation*}
a(u,u) \geq \alpha \left\| u \right\| ^{2}_{\mathcal{H}} \hspace{0.25in} \forall v \in \mathcal{H} 
\end{equation*}
The existence of the solution of the weak problem (\ref{weak_formulation}) and its uniqueness are the consequences of the Lax-Milgram Theorem \cite{Lax1954}, \cite{Ciarlet1991}. Note that the existence and uniqueness of the Galerkin weak problem (\ref{galerkin_weak_formulation}) is also a consequence of the Lax-Milgram theorem, since the bilinear form $a(\cdot,\cdot)$ restricted on $\mathcal{H}^{h}$ obviously remains bilinear, continuous and V-elliptic \cite{Solin2006}. Lax-Milgram theorem also states that the solution is bounded by the right hand side which represents the stability condition for the Galerkin method.

The convergence of the Galerkin approximation follows from C\'ea's lemma and an additional convergence theorem \cite{Cea1964}, \cite{Solin2006}. On the polygonal domains the sequence of subspaces of $\mathcal{H} = \mathbf{H}^{1}_{0} \left( B_{0}\right)^{9}$ can be obtained by the successive uniform refinement of the initial mesh using the midpoints as new nodes thus subdividing every triangle into 4 congruent triangles. Therefore $\mathcal{H}_{n} \subset \mathcal{H}_{n+1}$ for every $n \in \mathbb{N}$ and the sequence of spaces $\mathcal{H}^{n}$ is dense in $\mathcal{H}$ \cite{Ainsworth2000}, and thus
\begin{equation*}
\overline{\bigcup_{n=1}^{\infty} {\mathcal{H}_{n}}} = \mathcal{H}
\end{equation*}
and $u^{n}$ converges to $u$ as $n \rightarrow \infty$ \cite{Solin2006}, \cite{Duran2010}.

It was shown that there exists a sequence of triangulations that ensures optimal rates of convergence in $\mathbf{H}^{1}$-norm for the FEM approximation of the second order strongly elliptic system with zero Dirichlet boundary condition on polyhedron domain with continuous, piecewise polynomials of degree $m$ \cite{Bacuta2005}.

\vspace{.25in}

\section{Validation of the FEM for Different Boundary Conditions}

\vspace{.25in}

Let us consider the plate $B_{0}$ to be a square plate of size $[0,a]\times \lbrack 0,a]$ with the boundary $G=G_{1} \cup G_{2} \cup G_{3} \cup G_{4}$ and the hard simply supported boundary conditions written in terms of the kinematic variables in the mixed Dirichlet-Neumann:
\begin{eqnarray}
G_{1}\cup G_{2} &:&W=0,\text{ }W^{\ast }=0,\text{ }\Psi _{2}=0,\text{ }\Omega_{1}^{0}=0,\text{ }\hat{\Omega}_{1}^{0}=0,
\notag \\
&&\frac{\partial \Omega_{3}}{\partial n}=0,\text{ }\frac{\partial \Psi _{1}}{\partial n}=0,\text{ }\frac{\partial \Omega_{2}^{0}}{\partial n}=0,\text{ }\frac{\partial \hat{\Omega}_{2}^{0}}{\partial n}=0;  
\notag \\
G_{3}\cup G_{4} &:&W=0,\text{ }W^{\ast }=0,\text{ }\Psi _{1}=0,\text{ }\Omega_{2}^{0}=0,\text{ },\hat{\Omega}_{2}^{0}=0, \notag \\
&& \frac{\partial \Omega_{3}}{\partial n}=0,\text{ }\frac{\partial \Psi _{2}}{\partial n}=0,\text{ }\frac{\partial \Omega_{1}^{0}}{\partial n}=0,\text{ }\frac{\partial \hat{\Omega}_{1}^{0}}{\partial n}=0.  
\notag
\end{eqnarray}
where
\begin{eqnarray*}
G_{1} &=& \left\{ \left( x_{1},x_{2} \right) :x_{1} = 0 , x_{2}\in \left[ 0,a\right] \right\}, \notag \\ 
G_{2} &=& \left\{ \left( x_{1},x_{2} \right) :x_{1} = a , x_{2}\in \left[ 0,a\right] \right\}, \notag \\ 
G_{3} &=& \left\{ \left( x_{1},x_{2} \right) :x_{1}\in \left[ 0,a\right] , x_{2} = 0 \right\}, \notag \\
G_{4} &=& \left\{ \left( x_{1},x_{2} \right) :x_{1}\in \left[ 0,a\right] , x_{2} = a \right\}, \notag
\end{eqnarray*}

The existence of a sequence of triangulations that ensures the optimal rates of convergence for the Finite Element approximation of the solution of a second order strongly elliptic system with homogeneous Dirichlet boundary condition on polyhedron domain with continuous piecewise polynomials was shown in \cite{Bacuta2005}. For the case of piecewise linear polynomials the optimal rate of convergence in $\mathbf{H}^{1}$-norm is linear.

We propose to use the uniform refinement to form the sequence of triangulations and estimate the order of the error of approximation of the proposed FEM in $\mathbf{H}^{1}$-norm and $L_{2}$-norm.

Let us consider homogeneous Dirichlet boundary conditions. We will assume the solution $u$ of the form:
\begin{equation}
u_{i} = U_{i}\sin \left( \frac{\pi x_{1}}{a}\right) \sin \left( \frac{\pi x_{2}}{a}\right), \hspace{0.25in} U_{i} \in \mathbb{R}, i=\overline{1,9}, 
\label{test_solution}
\end{equation}
which automatically satisfies homogeneous Dirichlet boundary conditions. Substituting the solution (\ref{test_solution}) into the system of field equations (\ref{field_equations}) we can find the corresponding right-hand side function $f$. The results of the error estimation of the FEM approximation in $\mathbf{H}^{1}$ and $L_{2}$ norms performed for the elastic parameters corresponding to the polyurethane foam are given in the Tables \ref{tab:FEM_H1_Error_Dirichlet} and \ref{tab:FEM_L2_Error_Dirichlet} respectively.

Let us consider mixed Neumann-Dirichlet boundary conditions. Simply supported boundary conditions represent this type of boundary conditions and therefore the FEM approximation can be compared with the analytical solution developed in the Chapter 3 for some fixed value of the parameter $\eta$. The results of the error estimation of the FEM approximation in $\mathbf{H}^{1}$ and $L_{2}$ norms performed for the elastic parameters corresponding to the polyurethane foam are given in the Tables \ref{tab:FEM_H1_Error_Mixed} and \ref{tab:FEM_L2_Error_Mixed} respectively.
\noindent
\begin{table}
\caption {Order of Convergence in $\mathbf{H}^{1}$-norm for Homogeneous Dirichlet BC}
\label{tab:FEM_H1_Error_Dirichlet}
\setkeys{Gin}{keepaspectratio}
\resizebox*{\textwidth}{\textheight}
{\begin{tabular}{ccccc}
\\ [2ex]
\hline \\ 
Refinements & Number of Nodes & Diameter & Error in $\mathbf{H}^{1}$-norm & Convergence Rate \\ [2ex] \hline
\\ [0.5ex]
0 & 177 & 0.302456 & 1.620369 &  \\
[1ex]
1 & 663 & 0.151228 & 0.711098 & \bf{1.19} \\
[1ex]
2 & 2565 & 0.075614 & 0.322016 & \bf{1.14} \\
[1ex]
3 & 10089 & 0.037807 & 0.150149 & \bf{1.10} \\ 
[1ex]
4 & 40017 & 0.018903 & 0.073481 & \bf{1.03} \\
[1ex]
5 & 159393 & 0.009451 & 0.036512 & \bf{1.01} \\ [2ex]
\hline
\end{tabular}}
\end{table}

\noindent
\begin{table}
\caption {Order of Convergence in $L_{2}$-norm for Homogeneous Dirichlet BC}
\label{tab:FEM_L2_Error_Dirichlet}
\setkeys{Gin}{keepaspectratio}
\resizebox*{\textwidth}{\textheight}
{\begin{tabular}{ccccc}
\\ [2ex]
\hline \\ 
Refinements & Number of Nodes & Diameter & Error in $L_{2}$-norm & Convergence Rate \\ [2ex] \hline
\\ [0.5ex]
0 & 177 & 0.302456 & 0.279484 &  \\
[1ex]
1 & 663 & 0.151228 & 0.069632 & \bf{2.00} \\ 
[1ex]
2 & 2565 & 0.075614 & 0.018175 & \bf{1.94} \\ 
[1ex]
3 & 10089 & 0.037807 & 0.004598 & \bf{1.98} \\ 
[1ex]
4 & 40017 & 0.018903 & 0.001153 & \bf{2.00} \\
[1ex]
5 & 159393 & 0.009451 & 0.000288 & \bf{2.00} \\ [2ex] 
\hline
\end{tabular}}
\end{table}

\noindent
\begin{table}
\caption {Order of Convergence in $\mathbf{H}^{1}$-norm for Mixed Neumann-Dirichlet BC}
\label{tab:FEM_H1_Error_Mixed}
\setkeys{Gin}{keepaspectratio}
\resizebox*{\textwidth}{\textheight}
{\begin{tabular}{ccccc}
\\ [2ex]
\hline \\ 
Refinements & Number of Nodes & Diameter & Error in $\mathbf{H}^{1}$-norm & Convergence Rate \\ [2ex] \hline
\\ [0.5ex]
0 & 177 & 0.302456 & 0.236791 &  \\
[1ex]
1 & 663 & 0.151228 & 0.115809 & \bf{1.03} \\
[1ex]
2 & 2565 & 0.075614 & 0.054195 & \bf{1.09} \\
[1ex]
3 & 10089 & 0.037807 & 0.026233 & \bf{1.05} \\ 
[1ex]
4 & 40017 & 0.018903 & 0.012986 & \bf{1.01} \\
[1ex]
5 & 159393 & 0.009451 & 0.006475 & \bf{1.00} \\ [2ex]
\hline
\end{tabular}}
\end{table}

\noindent
\begin{table}
\caption {Order of Convergence in $L_{2}$-norm for Mixed Neumann-Dirichlet BC}
\label{tab:FEM_L2_Error_Mixed}
\setkeys{Gin}{keepaspectratio}
\resizebox*{\textwidth}{\textheight}
{\begin{tabular}{ccccc}
\\ [2ex]
\hline \\ 
Refinements & Number of Nodes & Diameter & Error in $L_{2}$-norm & Convergence Rate \\ [2ex] \hline
\\ [0.5ex]
0 & 177 & 0.302456 & $6.214 \times 10^{-2}$ &  \\
[1ex]
1 & 663 & 0.151228 & $1.638 \times 10^{-2}$ & \bf{1.92} \\ 
[1ex]
2 & 2565 & 0.075614 & $4.219 \times 10^{-3}$ & \bf{1.96} \\ 
[1ex]
3 & 10089 & 0.037807 & $1.065 \times 10^{-3}$ & \bf{1.99} \\ 
[1ex]
4 & 40017 & 0.018903 & $2.678 \times 10^{-4}$ & \bf{1.99} \\
[1ex]
5 & 159393 & 0.009451 & $6.772 \times 10^{-5}$ & \bf{1.98} \\ [2ex] 
\hline
\end{tabular}}
\end{table}

\vspace{.25in}

\subsection{Validation of the proposed FEM for Simply Supported Cosserat Elastic Plate}

\vspace{.25in}

The boundary condition for the variable $\Omega _{3}$ is a Neumann-type boundary condition:
\begin{equation*}
\frac{\partial \Omega_{3}}{\partial n}=0 \text{ on } G 
\end{equation*}
and thus we will look for $\Omega _{3}$ in the space $\mathbf{H}^{1}\left( \Delta, B_{0}\right)$, where. 
\begin{equation*}
\mathbf{H}^{1}\left( \Delta, B_{0}\right) = \left\{ u \in \mathbf{H}^{1} \left( B_{0}\right): \Delta u \in L_{2} (B_{0})\right\}
\end{equation*}

The boundary condition for the variables $W$ and $W^{\ast}$ is a Dirichlet-type boundary condition:
\begin{eqnarray*}
W&=&0 \text{ on } G \notag \\ 
W^{\ast}&=&0 \text{ on } G  \notag 
\end{eqnarray*}
and thus we will look for $W$ and $W^{\ast}$ in the space $\mathbf{H}^{1}_{0}\left( B_{0}\right)$ defined as  \cite{Johnson1987}:
\begin{equation*}
\mathbf{H}^{1}_{0}\left( B_{0}\right) = \left\{v\in \mathbf{H}^{1}\left( B_{0}\right) ,v=0 \text{ on } G \right\} \end{equation*}

The boundary condition for the variables $\Psi _{1}$, $\Omega_{2}^{0}$ and $\hat{\Omega}_{2}^{0}$ is of mixed Dirichlet-Neumann type:
\begin{eqnarray*}
\frac{\partial \Psi _{1}}{\partial n}=0, \frac{\partial \Omega_{2}^{0}}{\partial n}=0, \frac{\partial \hat{\Omega}_{2}^{0}}{\partial n}=0 &&\text{ on } G_{1}\cup G_{2}\notag \\ 
\Psi _{1}=0, \Omega_{2}^{0}=0, \hat{\Omega}_{2}^{0}=0 &&\text{ on } G_{3}\cup G_{4}  \notag 
\end{eqnarray*}
and thus we will look for $\Psi _{1}$, $\Omega_{2}^{0}$ and $\hat{\Omega}_{2}^{0}$ in the following space \cite{Johnson1987}:
\begin{equation*}
\mathbf{H}^{1}_{V} = \left\{ v\in \mathbf{H}^{1}\left( \Delta, B_{0}\right) , v=0 \text{ on } G_{3}\cup G_{4} \right\}
\end{equation*}

The boundary condition for the variables $\Psi_{2}$, $\Omega_{1}^{0}$ and $\hat{\Omega}_{1}^{0}$ is of mixed Dirichlet-Neumann type:
\begin{eqnarray*}
\Psi _{2}=0, \Omega_{1}^{0}=0, \hat{\Omega}_{1}^{0}=0 &&\text{ on } G_{1}\cup G_{2} \notag \\ 
\frac{\partial \Psi _{2}}{\partial n}=0, \frac{\partial \Omega_{1}^{0}}{\partial n}=0, \frac{\partial \hat{\Omega}_{1}^{0}}{\partial n}=0 &&\text{ on } G_{3}\cup G_{4} \notag 
\end{eqnarray*}
and thus we will look for $\Psi _{2}$, $\Omega_{1}^{0}$ and $\hat{\Omega}_{1}^{0}$ in the following space \cite{Johnson1987}:
\begin{equation*}
\mathbf{H}^{1}_{H} = \left\{ v\in \mathbf{H}^{1}\left( \Delta, B_{0}\right) , v=0 \text{ on } G_{1}\cup G_{2} \right\}
\end{equation*}

Therefore we will look for the solution
\begin{equation*}
\left[ \Psi _{1},\Psi _{2},W,\Omega_{3},\Omega _{1}^{0},\Omega _{2}^{0},W^{\ast },\hat{\Omega}_{1},\hat{\Omega}_{2}\right] ^{T}
\end{equation*}
of the Cosserat plate field equations (\ref{field_equations}) in the space $\mathcal{H}$ defined as
\begin{equation}
\mathcal{H} = \mathcal{H}_{1} \times \mathcal{H}_{2} \times \mathcal{H}_{3} \times \mathcal{H}_{4} \times \mathcal{H}_{5} \times \mathcal{H}_{6} \times \mathcal{H}_{7} \times \mathcal{H}_{8} \times \mathcal{H}_{9}
\end{equation}
where 
\begin{eqnarray}
&& \mathcal{H}_{1} = \mathcal{H}_{6} = \mathcal{H}_{9} = \mathbf{H}^{1}_{V}\left( B_{0}\right), \notag \\
&& \mathcal{H}_{2} = \mathcal{H}_{5} = \mathcal{H}_{8} = \mathbf{H}^{1}_{H}\left( B_{0}\right), \notag \\
&& \mathcal{H}_{3} = \mathcal{H}_{7} = \mathbf{H}^{1}_{0}\left( B_{0}\right), \notag \\ 
&& \mathcal{H}_{4} = \mathbf{H}^{1}\left( \Delta, B_{0}\right). \notag 
\end{eqnarray}

The space $\mathcal{H}$ is a Hilbert space equipped with the inner product $\left\langle u,v \right\rangle_{\mathcal{H}}$ on defined on $\mathcal{H}$ as follows:
\begin{equation*}
\left\langle u,v \right\rangle_{\mathcal{H}} = \sum_{i=1}^{9}{\left\langle u_{i},v_{i} \right\rangle_{\mathcal{H}_{i}}} \text{ for } u,v \in \mathcal{H}
\end{equation*}
where $\left\langle u,v \right\rangle_{\mathcal{H}_{i}}$ is an inner product defined on the Hilbert space $\mathcal{H}_{i}$ respectively. 

Taking into account the essential boundary conditions we define the finite element spaces $\mathcal{H}_{i}^{h}$ as follows:
\begin{eqnarray}
&& \mathcal{H}_{1}^{h} = \mathcal{H}_{6}^{h} = \mathcal{H}_{9}^{h} = \left\{ v: v \in C\left(B_{0}\right), v \text{ is linear on every } K_{j}, v=0 \text{ on } G_{1} \cup G_{2} \right\}, \notag \\
&& \mathcal{H}_{2}^{h} = \mathcal{H}_{5}^{h} = \mathcal{H}_{8}^{h} = \left\{ v: v \in C\left(B_{0}\right), v \text{ is linear on every } K_{j}, v=0 \text{ on } G_{3} \cup G_{4} \right\}, \notag \\
&& \mathcal{H}_{3}^{h} = \mathcal{H}_{7}^{h} = \left\{ v: v \in C\left(B_{0}\right), v \text{ is linear on every } K_{j}, v=0 \text{ on } G \right\}, \notag \\ 
&& \mathcal{H}_{4}^{h} = \left\{ v: v \in C\left(B_{0}\right), v \text{ is linear on every } K_{j} \right\}. \notag
\end{eqnarray}

The finite dimensional space $\mathcal{H}^{h}$ is then defined as
\begin{equation}
\mathcal{H}^{h} = \mathcal{H}_{1}^{h} \times \mathcal{H}_{2}^{h} \times \mathcal{H}_{3}^{h} \times \mathcal{H}_{4}^{h} \times \mathcal{H}_{5}^{h} \times \mathcal{H}_{6}^{h} \times \mathcal{H}_{7}^{h} \times \mathcal{H}_{8}^{h} \times \mathcal{H}_{9}^{h}
\end{equation}

We solve the field equations using described Finite Element method and compare the obtained results with the analytical solution for the square plate made of polyurethane foam derived in the Chapter 3.
\noindent
\begin{table}
\caption {Order of Convergence in $\mathbf{H}^{1}$-norm for Simply Supported Plate}
\label{tab:FEM_H1_Error_Supported}
\setkeys{Gin}{keepaspectratio}
\resizebox*{\textwidth}{\textheight}
{\begin{tabular}{ccccc}
\\ [2ex]
\hline \\ 
Refinements & Nodes Number & Diameter & Error in $\mathbf{H}^{1}$-norm & Convergence Rate \\ [2ex] \hline
\\ [0.5ex]
0 & 177 & 0.302456 & 0.256965 &  \\ 
[1ex]
1 & 663 & 0.151228 & 0.119234 & \bf{1.11} \\ 
[1ex]
2 & 2565 & 0.075614 & 0.054701 & \bf{1.12} \\
[1ex]
3 & 10089 & 0.037807 & 0.026301 & \bf{1.05} \\ 
[1ex]
4 & 40017 & 0.018903 & 0.012994 & \bf{1.01} \\ 
[1ex]
5 & 159393 & 0.009451 & 0.006476 & \bf{1.00} \\ [2ex] 
\hline
\end{tabular}}
\end{table}
\noindent
\begin{table}
\caption {Order of Convergence in $L_{2}$-norm for Simply Supported Plate}
\label{tab:FEM_L2_Error_Supported}
\setkeys{Gin}{keepaspectratio}
\resizebox*{\textwidth}{\textheight}
{\begin{tabular}{ccccc}
\\ [2ex]
\hline \\ 
Refinements & Nodes Number & Diameter & Error in $L_{2}$-norm & Convergence Rate \\ [2ex] \hline
\\ [0.5ex]
0 & 177 & 0.302456 & $8.253 \times 10^{-2}$ &  \\ 
[1ex]
1 & 663 & 0.151228 & $2.260 \times 10^{-2}$ & \bf{1.87} \\
[1ex]
2 & 2565 & 0.075614 & $5.860 \times 10^{-3}$ & \bf{1.95} \\ 
[1ex]
3 & 10089 & 0.037807 & $1.482 \times 10^{-3}$ & \bf{1.98} \\ 
[1ex]
4 & 40017 & 0.018903 & $3.720 \times 10^{-4}$ & \bf{1.99} \\ 
[1ex]
5 & 159393 & 0.009451 & $9.355 \times 10^{-5}$ & \bf{1.99} \\ [2ex]
\hline
\end{tabular}}
\end{table}

The initial distribution of the pressure, as in the Chapter 4, is assumed sinusoidal:
\begin{equation}
p\left( x_{1},x_{2}\right) =\sin \left( \frac{\pi x_{1}}{a}\right) \sin
\left( \frac{\pi x_{2}}{a}\right)  \label{initial_sinusoidal_pressure}
\end{equation}

The estimation of the error in $\mathbf{H}^{1}$ norms shows that the order of the error is optimal (linear) in $\mathbf{H}^{1}$-norm for the piecewise linear elements for the simply supported plate. The results of the error estimation of the FEM approximation in $\mathbf{H}^{1}$ and $L_{2}$ norms performed for the elastic parameters corresponding to the polyurethane foam are given in the Tables \ref{tab:FEM_H1_Error_Supported} and \ref{tab:FEM_L2_Error_Supported} respectively.

The comparison of the maximum of the displacements $u_{i}$ and microrotations $\varphi_{i}$ calculated using Finite Element method with 320 thousand elements and the analytical solution for the micropolar plate theory is provided in the Table \ref{tab:Relative_Error}. The relative error of the approximation of the optimal value of the splitting parameter is $0.09\%$.
\noindent
\begin{table}
\caption {Relative Error of the Maximum Values of the Displacement and Microrotations}
\label{tab:Relative_Error}
\setkeys{Gin}{keepaspectratio}
\resizebox*{\textwidth}{\textheight}
{\begin{tabular}{ccccccc}
\\ [2ex]
\hline \\ 
& Optimal $\eta$ & $u_{1}$ & $u_{2}$ & $u_{3}$ & $\varphi_{1}$ & $\varphi_{2}$ \\ [2ex] \hline
\\ [0.5ex]
Finite Element Solution & 0.040760 & -0.014891 &  -0.014891 & 0.307641 & 0.046767 &  -0.046767 \\
[1ex]
Analytical Solution & 0.040799 & -0.014892 & -0.014892 & 0.307674 & 0.046770 & -0.046770 \\
[1ex]
Relative Error (\%) & \bf{0.09} & \bf{0.03} & \bf{0.03} & \bf{0.04} & \bf{0.03} & \bf{0.03} \\ [2ex] 
\hline
\end{tabular}}
\end{table}

\begin{figure}
\begin{center}
	\includegraphics[width=5in]{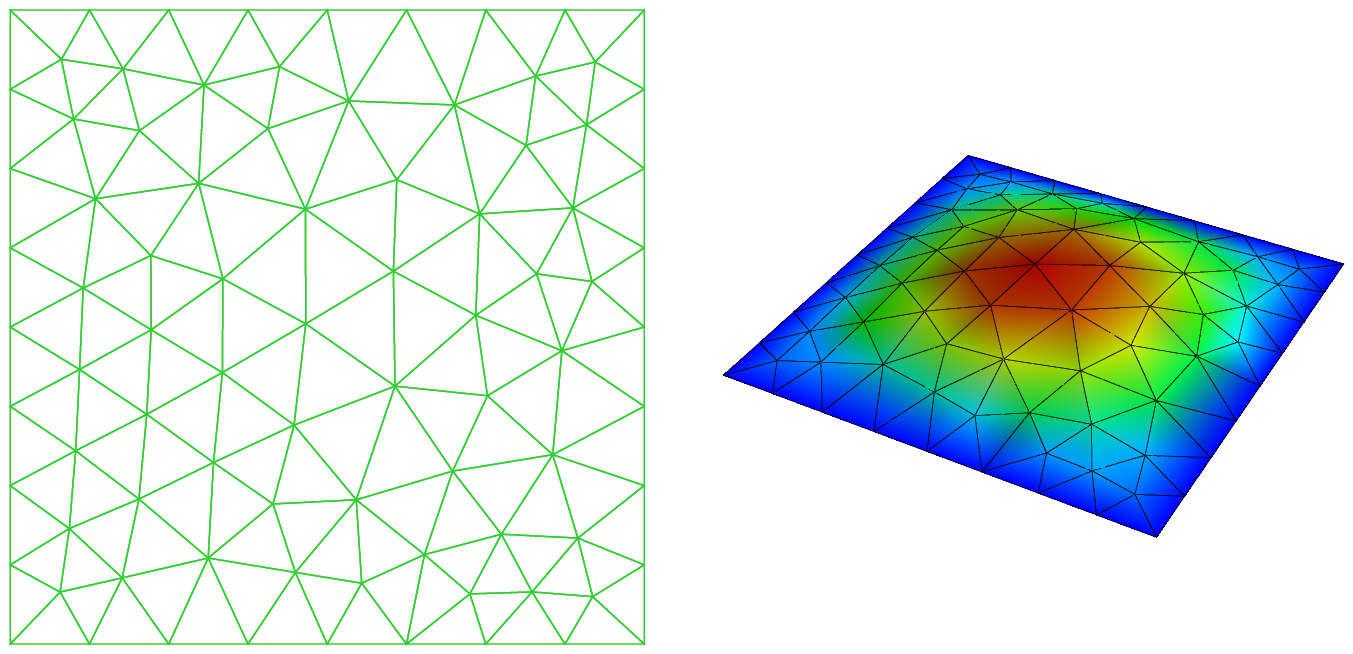}
	\vspace{.1in}
	\caption{Hard simply supported square plate 2.0m$\times $2.0m$\times $0.1m made of polyurethane foam: the initial mesh and the isometric view of the resulting vertical deflection of the plate}
	\label{fig:Plate_0}
\vspace{.15in}
\end{center}
\end{figure}

\begin{figure}
\begin{center}
	\includegraphics[width=6in]{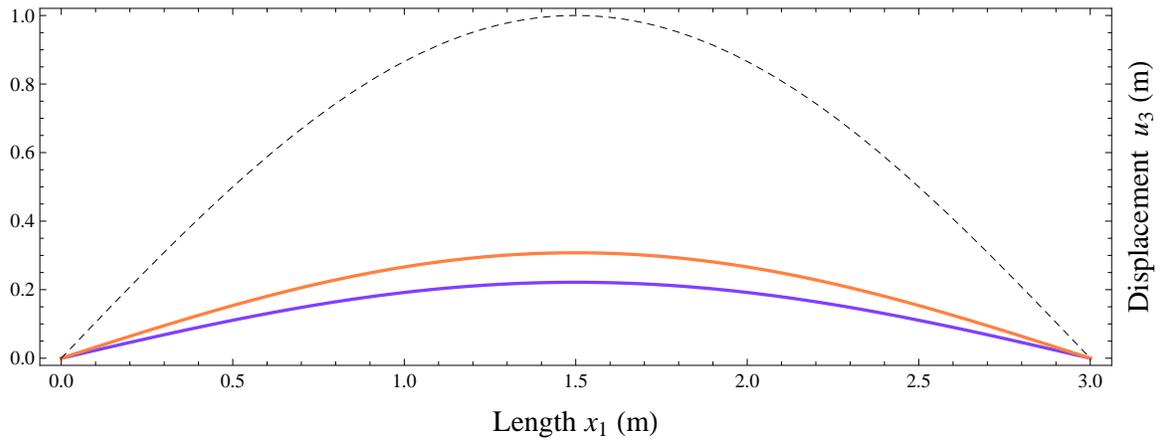}
	\vspace{.1in}
	\caption{The cross section that contains the center of the micropolar square plate 3.0m$\times$3.0m$\times$0.1m made of polyurethane foam under the sinusoidal load: Cosserat clamped plate -- solid blue line, Cosserat simply supported plate -- solid orange line, initial sinusoidal load -- dashed black line.}
	\label{fig:plates_cross_comparison}
\end{center}
\end{figure}

The Figure \ref{fig:Plate_0} represents the Finite Element modeling of the bending of the simply supported square plate made of polyurethane foam. The comparison of the distribution of the vertical deflection of the clamped and simply supported plates is given in the Figure \ref{fig:plates_cross_comparison}.

\noindent

\noindent
\noindent
\begin{figure}
\begin{center}
	\includegraphics[width=6in]{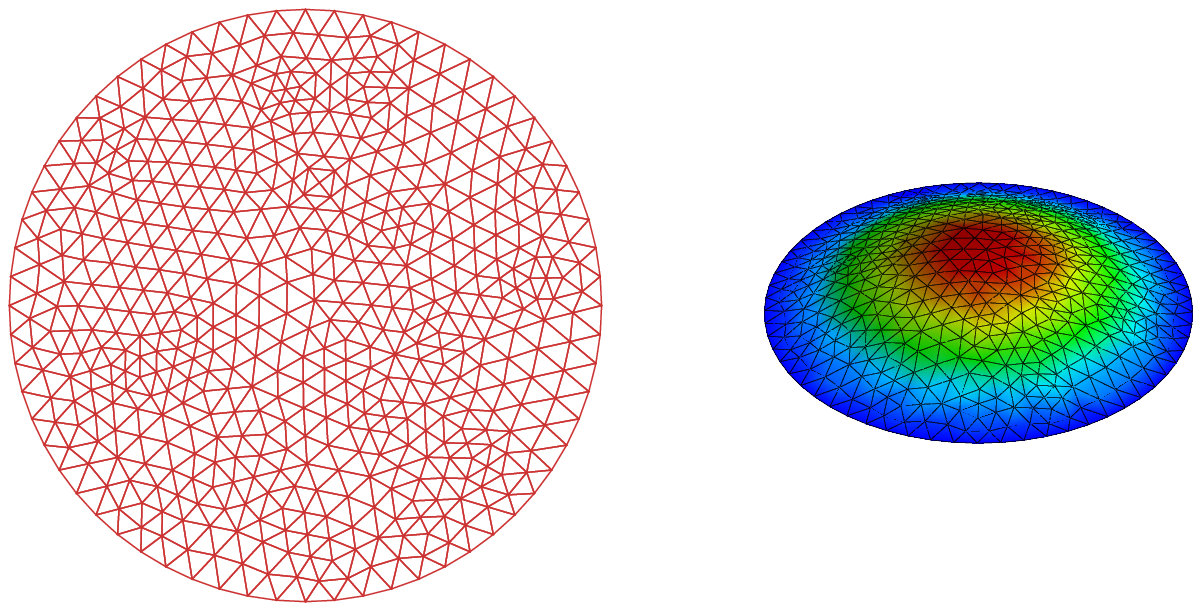}
	\vspace{.1in}
	\caption{The circular Cosserat clamped plate of radius $R=1.0$m and thickness $h=0.1$m made of polyurethane foam under the uniform load: the initial mesh and the isometric view of the resulting vertical deflection of the plate.}
	\label{fig:circular_plate_1}
\end{center}
\end{figure}
\noindent
\begin{figure}
\begin{center}
	\includegraphics[width=6in]{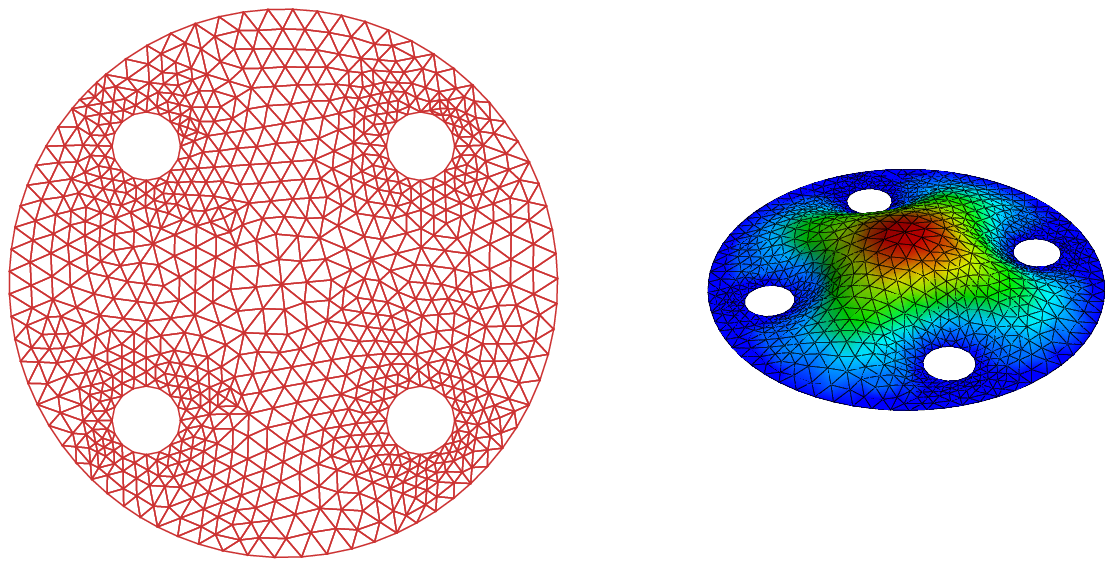}
	\vspace{.1in}
	\caption{The circular Cosserat clamped plate of radius $R=1.0$m and thickness $h=0.1$m made of polyurethane foam with circular clamped holes under the uniform load: the initial mesh and the isometric view of the resulting vertical deflection of the plate.}
	\label{fig:circular_plate_2}
\end{center}
\end{figure}
\noindent
\begin{figure}
\begin{center}
	\includegraphics[width=6in]{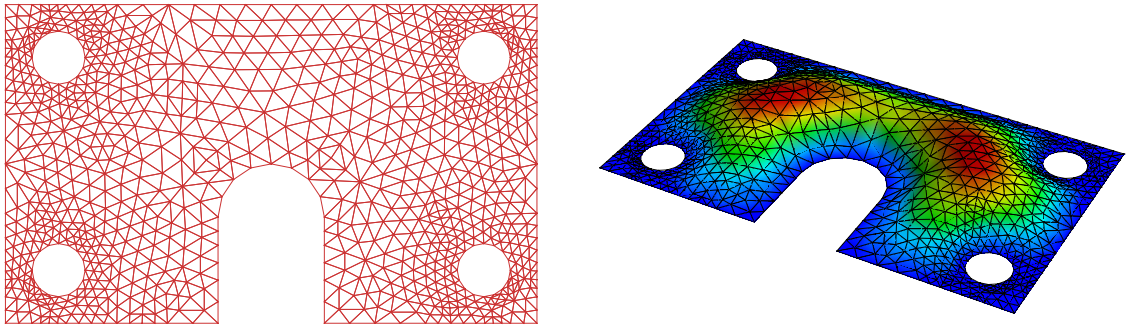}
	\vspace{.1in}
	\caption{The clamped plate of size 10.0m$\times$6.0m$\times$0.1m made of polyurethane foam under the uniform load: the initial mesh and the isometric view of the resulting vertical deflection of the plate.}
	\label{fig:weird_plate_1}
\end{center}
\end{figure}
\noindent
\begin{figure}
\begin{center}
	\includegraphics[width=6in]{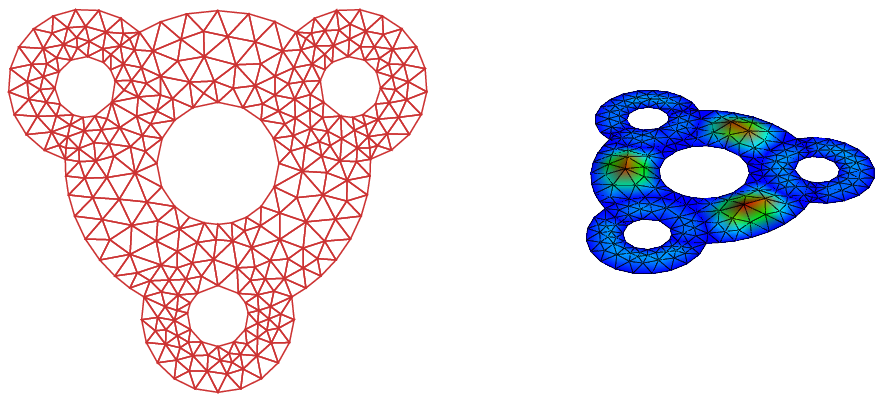}
	\vspace{.1in}
	\caption{The clamped polyurethane gasket under the uniform load: the initial mesh and the isometric view of the resulting vertical deflection of the plate.}
	\label{fig:weird_plate_3}
\end{center}
\end{figure}

\vspace{.25in}

\section{Conclusion}

This article develops and validates the Finite Element Method for Cosserat elastic plates based on the enhanced Cosserat plate theory. We present the Finite Element analysis of the Cosserat plates of different shapes, under different loads and different boundary conditions. We discuss the existence and uniqueness of the weak solution and the convergence of the proposed FEM. The proposed finite element method yields an optimum convergence rate, when comparing the main kinematic variables with the analytical solution of the two-dimensional problem. We also consider the numerical analysis of plates with circular holes. We calculate the stress concentration factor around the hole and show that it is smaller would be expected on the basis of Reissner theory for simple elastic plates. The finite element comparison of the plates with holes confirm that smaller holes exhibit less stress concentration.

\end{document}